\documentclass[12pt]{article}
\usepackage[cp1251]{inputenc}
\usepackage[T2A]{fontenc}
\usepackage[english]{babel}

\def\eps{\varepsilon}

\newcounter{num}[section]

\newcommand{\Th}{\refstepcounter{num}
{\bf Theorem \arabic{section}.\arabic{num} }}

\newcommand{\Lemma}{\refstepcounter{num}
{\bf Lemma \arabic{section}.\arabic{num} }}

\newcommand{\Cor}{\refstepcounter{num}
{\bf Corollary \arabic{section}.\arabic{num} }}

\newcommand{\Note}{\refstepcounter{num}
{\it Note \arabic{section}.\arabic{num} }}

\newcommand{\St}{\refstepcounter{num}
{\bf Statement \arabic{section}.\arabic{num} }}

\newcommand{\Def}{\refstepcounter{num}
{\it Definition \arabic{section}.\arabic{num} }}

\newcommand{\Proof}{{\bf Proof. }}

\def\v{\vec}

\def\a{\alpha}
\def\d{\delta}
\def\m{\times}

\def\t{\tilde}

\def\Z_N{{\bf Z}_N}

\def\_phi{\varphi}

\def\diam{{\rm diam}}

\author{Shkredov I.D.}
\title{On Dynamical Systems With Slow Recurrence Time
\footnote{This work was supported by the program
"Leading Scientific Schools" (project no. 136.2003.1),
by RFFI grant no. 02--01--00912
and
INTAS (grant no. 03--51--5-70).}.}
\date{}
\begin{document}
\maketitle

\refstepcounter{section}

{\bf \arabic{section}. Introduction.}

Let  {\it Х} be a set with Borel sigma--algebra of measurable sets $\mathcal{B}$.
Let $T$ be a transformation of $X$ into itself,
preserving a measure $\mu$
and let us assume that measure of $X$ is equal to $1$.
$(X, \mathcal{B}, \mu, T)$ is called a measure--preserving dynamical system.
The well--known
{\bf H. Poincar\'{e}'s  Recurrence Theorem }
\cite{Poin}
asserts
that for any measurable set $E \subseteq X$,
$\mu E >0$
there exists a natural number $n$ such that $\mu(E \cap T^{-n} E) > 0$.

Suppose, in addition, $X$ is a metric space with metric
$ d(\cdot ,\cdot )$.
In the case Poincar\'{e}'s  Theorem can be reformulate  as

\Th {\bf (H. Poincar\'{e})}
{\it
  Let {\it X } be a metric space with metric $d(\cdot,\cdot)$
  and Borel sigma--algebra of measurable sets $\mathcal{B}$.
  Suppose ${\it T}$ is a measure--preserving transformation of $X$ into itself.
%  preserving a measure $\mu$.
  Then for almost every point $ x \in X $ the following inequality holds
    $$
        \liminf_{n \to \infty}~  d(T^n x,x)  = 0 \,.
    $$
}
\label{t:Poin}

Poincar\'{e}'s Theorem was generalized by  H. Furstenberg, Y. Katznelson and D. Ornstein in \cite{Fu,Fu_KO}.

\Th {\bf (H. Furstenberg, Y. Katznelson, D. Ornstein)}
{\it
 Let  $X$ be a set with Borel sigma--algebra of measurable sets $\mathcal{B}$.
 Let ${\it T}$ be a measure--preserving transformation of $X$ into itself
 % preserving a measure $\mu$
 and $k\ge 3$.
 Then for any measurable set $E$, $\mu E > 0$
 there exists a natural number $n$ such that
 $$
   \mu(E \cap T^{-n} E \cap T^{-2n} E  \cap \dots \cap T^{-(k-1) n} E) > 0 \,.
 $$
}
\label{t:F_deg}

If $X$ is a metric space then Theorem \ref{t:F_deg} can be rewrite as

{\bf Theorem \arabic{section}.\arabic{num}$'$} {\bf (H. Furstenberg, Y. Katznelson, D. Ornstein)}
{\it
Let $X$ be a metric space with metric
$d(\cdot,\cdot)$ and Borel sigma--algebra of measurable sets
$\mathcal{B}$.
Suppose ${\it T}$ is a measure--preserving transformation of $X$ into itself
%preserving a measure $\mu$
and $k\ge 3$.
Then for almost every point $x\in X$
%the following inequality holds
$$
  \liminf_{n\to \infty}~ \max
  \{ d(T^{n}x,x), d(T^{2n}x,x), \dots, d(T^{(k-1)n}x,x) \} = 0 \,.
$$
}
\label{t:F_deg_m}

In fact, Furstenberg proved a more general result.
He proved that powers of transformation  $T$ in Theorem \ref{t:F_deg}
can be replaced by any finite number of  commutative transformations.

\Th {\bf (H. Furstenberg)}
{\it
 Let $X$ be a metric space with metric
$d(\cdot,\cdot)$ and Borel sigma--algebra of measurable sets
$\mathcal{B}$.
Suppose  $k\ge 2$ and $T_1, T_2, \dots T_k$ are {\it commutative} measure--preserving
transformations.
%preserving a measure $\mu$.
Then for almost every point $x\in X$
%the following inequality holds
$$
  \liminf_{n\to \infty}~ \max
  \{ d(T_1^{n}x,x), d(T_2^{n}x,x), \dots, d(T_{k}^{n}x,x) \} = 0 \,.
$$
}
\label{t:F_comm}

Let $A$ be a subset of positive integers.
By $[N]$ denote the segment $\{ 1, 2, \dots, N \}$.
The {\it upper density} of a set $A$ is defined to be
$$
  D^* (A) = \limsup_{N\to \infty} \frac{| A \cap [N] |}{N} \,\,,
$$
where $|A \cap [N]|$ is the cardinality of $A \cap [N]$.

As was showed in \cite{Fu}  Theorem \ref{t:F_deg} is equivalent of the famous E. Szemer\'{e}di's Theorem
on arithmetic progressions.

\Th {\bf (E. Szemer\'{e}di)}
{\it
 Let $A$ be a subset of positive integers and $D^* (A) > 0$.
 Then for any natural number $k\ge 3$
 the set $A$ contains an arithmetic progression of the length $k$.
}
\label{t:Sz}

It is easy to see that Theorem \ref{t:Sz}
implies Theorems \ref{t:F_deg} and  \ref{t:F_deg_m}$'$ (see \cite{Fu}).
%In fact, Theorem \ref{t:F_deg} (\ref{t:F_deg_m}$'$) is equivalent of Theorem \ref{t:Sz}.
In fact, Theorem \ref{t:F_deg} (\ref{t:F_deg_m}$'$) and Theorem \ref{t:Sz} are equivalent statements.
To prove this Furstenberg obtained the following beautiful result
which is called {\it Furstenberg's Correspondence Principle}.

\Th {\bf (H. Furstenberg)}
{\it
  Let  $A$ be a subset of positive integers with $D^* (A) > 0$.
  Then there exists a dynamical system
  $(X, \mathcal{B}, \mu, T)$ and a measurable set $E$, $\mu E = D^* (A)$
  such that for all integers $k\ge 3$ and all integers
  $m_1,m_2,\dots,$                                                                         % \\
  $m_{k-1} \ge 1$,
  \begin{equation}\label{e:D>mu}
    D^* ( A \cap (A+m_1) \cap \dots \cap (A+m_{k-1}) ) \ge \mu (E \cap T^{-m_1} E \cap \dots \cap T^{-m_{k-1}} E) \,.
  \end{equation}
}
\label{t:F_corr}
Theorem \ref{t:F_corr}
shows that there is a close connection
between Ergodic Theory and some combinatorial problems concerning  arithmetic progressions.

\St {\bf (H. Furstenberg)}
{\it
  Theorems \ref{t:F_deg} and \ref{t:F_corr} imply Theorem \ref{t:Sz}.
}
\label{st:F_Sz_derive}
\\
\Proof
Let $k$ be a natural number, $k\ge 3$.
Suppose that  a set  $A \subseteq {\bf N}$ does not contain arithmetic progressions
of the length $k$ and has positive upper density.
Using Theorem \ref{t:F_corr},  we obtain a dynamical system
$(X, \mathcal{B}, \mu, T)$ and a measurable set $E$ of positive measure such that
for all integers $m_1,m_2,\dots, m_{k-1} \ge 1 $ the inequality (\ref{e:D>mu}) holds.
On the other hand, using Theorem \ref{t:F_deg}, we get a natural $n$ such that
\begin{equation}\label{e:tmp_3:15:15}
   \mu(E \cap T^{-n} E \cap T^{-2n} E  \cap \dots \cap T^{-(k-1) n} E) > 0 \,.
\end{equation}
Put $m_1 = n, m_2 = 2n, \dots, m_{k-1} = (k-1)n$.
Then
(\ref{e:D>mu}) and (\ref{e:tmp_3:15:15}) imply
%inequality
$D^* ( A \cap (A+n) \cap \dots \cap ( A+(k-1)n ) ) > 0$.
%But
By assumption
$A$ does not contain arithmetic progressions
of the length $k$.
This completes the proof.
\\

The main goal  of this article is to obtain quantitative analog of Theorem \ref{t:F_corr}.
To formulate our result we need in some definitions.

Let us consider
a measure $H_h (\cdot)$ on
$X$, defined as
\begin{equation}\label{f:Hausdorff}
  H_h (E)= \lim_{\delta \rightarrow 0} H_{h}^{\delta }(E) \,,
\end{equation}
where $h(t)$ is a positive ($ h(0)=0 $) continuous increasing
function and $H_{h}^{\delta }(E)= \inf_\tau \{ \sum h(\delta_j)
\}$, when $\tau$ runs through all countable coverings $E$ by open
sets
$ \{ B_j \}$ , $diam(B_j) = \delta_j < \delta$.\\
If $h(t) = t^\alpha$, then we get the ordinary Hausdorff measure
$H_{\alpha } (\cdot ) $.

Generally speaking,
$H_h (\cdot)$ is an outer measure but
it is sigma--additive measure on Carath\'{e}odory's sigma--algebra of measurable sets
(see \cite{Bo} for details).
It is well--known that this sigma--algebra contains all Borel sets.

We shall say that a measure $\mu$ is congruent to a measure $H_h$, if
any $\mu$--measurable set is $H_h$--measurable (in the sense of Carath\'{e}odory)

\Def
Let
$$
     C(x) = \liminf_{n \to \infty}  \{ n \cdot h(d(T^{n}x,x)) \} \,.
$$
The function $C(x)$ is called {\it constant of recurrence}  for point $x$.

The first quantitative version of Theorem \ref{t:Poin} was proved by M. Boshernitzan in \cite{Bo}.
(A similar result was obtained independently by N.G. Moshchevitin in \cite{Mo}).

\Th {\bf (M. Boshernitzan)}
{\it
    Let {\it X } be a metric space with $ H_{h} (X) < \infty $
    and  ${\it T}$ be a measure--preserving transformation of $X$ into itself.
    Assume that $\mu$ is congruent to $H_h$.
    Then for almost every point  $x \in X $
    we have
    $C(x) < \infty$.
}
\label{t:Bo}

The following result (see \cite{Sh}) improves
%slightly
Theorem \ref{t:Bo}.

\Th
{\it
    Let {\it X } be a metric space with $ H_{h} (X) < \infty $
    and  ${\it T}$ be a measure--preserving transformation of $X$ into itself.
    Assume that $\mu$ is congruent to $H_h$.
    Then the function $C(x)$ is $\mu$--integrable
    and for any $\mu$--measurable set $A$, we have
    \begin{equation}\label{e:integral}
         \int_A C(x) d\mu \le H_h (A).
    \end{equation}
    If $H_h(A)=0$, then $\int_A C(x) d\mu = 0$
    with no demand on measures $\mu$ and $H_h$ to be congruent.
}
\label{t:Sh1}

\Note
    According to an example from \S 7 of paper \cite{Bo} the inequality (\ref{e:integral})
    %in Theorem \ref{t:Sh1}
    is best possible.

Let us return to Theorems \ref{t:F_deg} and \ref{t:Sz}.

 Suppose $N$ and $k$ are natural numbers, $k\ge 3$.
 We set
$$
  a_k(N) = \frac{1}{N} \max \{ |A| ~:~ A \subseteq [1,N],
$$
$$
  A \mbox{ contains no arithmetic progressions of length } k
   \},
$$
Clearly, Theorem \ref{t:Sz} can be reformulate as
$\lim_{N\to \infty} a_k (N) = 0$ for all $k\ge 3$.

The first result concerning the rate at which $a_k(N)$ approaches zero in the case of $k=3$
was obtained by K.F. Roth
(see \cite{Rt}).
In his paper Roth used the Hardy -- Littlewood method
to prove the inequality
$$
  a_3(N) \ll \frac{1}{\log \log N}.
$$
In other words Roth obtained a quantitative version of Theorem
\ref{t:Sz} and therefore Theorem \ref{t:F_deg} in the case of $k$ equals three.

   At present, the best upper bound for $a_3(N)$ is due to J. Bourgain.
He proved that
\begin{equation}
 a_3(N) \ll \sqrt{ \frac{\log \log N}{ \log N} } \,.
\end{equation}

In \cite{Gow_m} W.T. Gowers
obtained a quantitative result concerning the rate at which $a_k(N)$ approaches zero for all $k\ge 4$.

\Th {\bf (W.T. Gowers)}
{\it
  For all $k\ge 4$ the following inequality holds
  $a_k(N) \ll 1/ (\log \log N)^{c_k}$,
  where $c_k$ is an absolute constant depends on $k$ only.
}

A. Behrend \cite{Be} obtained a lower bound for $a_3(N)$.
His result was generalized by R. Rankin in \cite{Ra}
in the case of all $k\ge 3$ (see also \cite{Laba}).

\Th {\bf (A. Behrend, R. Rankin)}
{\it
Let $\eps>0$ be an arbitrary real number and $k\ge 3$ be a natural.
Then for all sufficiently large $N$ the following inequality holds
$$
  a_k (N) \ge \exp( - (1+\eps) C_k (\log N)^{1/(k+1)} ) \,,
$$
where $C_k$ is an absolute positive constant depends on $k$ only.
}
\label{t:Be_Ra}

In the case of $k=2$ quantitative version of Theorem \ref{t:F_comm} was obtained
in \cite{Shkr_tri_DAN,Shkr_tri} and improved in \cite{Shkr_tri_B}.
  Consider the two--dimensional lattice $[1,N]^2$ with basis  $\{(1,0)$, $(0,1)\}$.
  Let
$$
  L(N) = \frac{1}{N^2} \max \{~ |A| ~:~ A\subseteq [1,N]^2 ~\mbox{ and }~
$$
\begin{equation}\label{tri}
  A \mbox { contains no triples of the form } (k,m),~ (k+d,m),~ (k,m+d),~ d>0 \}.
\end{equation}

\Th
{\it
%  The following inequality holds
  We have
  $L(N) \ll 1/ (\log \log N)^{C'}$,
  where $C'$ is an absolute constant.
}

This inequality implies a result concerning recurrence time in Theorem \ref{t:F_comm}
in the case of $k=3$ (see \cite{Shkr_mult}).
Suppose
$S$ and $R$ are  two {\it commutative} measure--preserving transformations of $X$.

\Def
By $C_{S,R}(x)$ denote the function
$$
 C_{S,R}(x) =
    \liminf_{n \to \infty}~  \{L^{-1}(n) \cdot \max \{ h(d(S^{n}x,x)), h(d(R^{n}x,x)) \} \} \,,
$$
where $L^{-1}(n) = 1 / L(n)$.
$C_{S,R}(x)$ is called {\it constant of multiple
recurrence
%одновременного (или кратного)
} for point $x$.
\label{const_m_rec}

\Th
{\it Let {\it X } be a metric space with
$ H_h (X) = C < \infty $
and let  $S,R$
be two  commutative measure--preserving transformation of $X$.
Assume that $\mu$ is congruent to $H_h$.
Then the function $C_{S,R}(x)$ is $\mu$--integrable
and for any $\mu$--measurable set $A$
$$
  \int_A C_{S,R}(x) d\mu \le H_h (A).
$$
If $H_h(A)=0$ then $\int_A C_{S,R}(x) d\mu = 0$
with no demand on measures $\mu$ and $H_h$ to be congruent.
}
\label{t:int_2}

Let us formulate our main result.

Let $k$ be a natural number, $k\ge 3$.
Suppose that for any natural $N$ there exists a nonempty set
$A^{(N)} \subseteq \Z_N$ without arithmetic progressions
of the length $k$.
By $\rho(N)$ denote the density of  $A^{(N)}$ in $\Z_N$, $\rho(N) = |A^{(N)}| / N$.
We have $\rho(1)=1$.
Since $A^{(N)}$ contains no arithmetic progressions in $\Z_N$
it follows that
$A^{(N)}$ has no arithmetic progressions in ${\bf Z}$.
By Theorem \ref{t:Sz} we have $\rho(N) \to 0$ as $N \to \infty$.
Assume that $\rho(N)$ is a non--increasing  function.

\Th
{\it
  Let $\psi : {\bf N} \to {\bf R^+}$ be a monotonically increasing
  function, $X = [0,1]$ and $\mu$ is Lebesgue measure on $X$.
  Then there exists a dynamical system $(X,\mathcal{B},\mu,T,d)$ such that
  $\mu$ is congruent to Hausdorff measure $H_1$,
  $H_1 (X) = 0$
  and
  for almost every point $x\in X$
  \begin{equation}\label{i:main_rec}
    \liminf_{n\to \infty}~ \{ \rho^{-1} (n) \psi (n) \max  \{ d(T^n x,x), d(T^{2n} x,x), \dots, d(T^{(k-1)n} x,x) \} \} \ge 1 \,,
  \end{equation}
  where
  $
   \rho^{-1} (n) = 1 / \rho (n).
  $
}
\label{t:general_case}

\Note The identity  $H_1 (X) = 0$ in Theorem \ref{t:general_case}
is very important.
If this equality is not true then our Theorem is trivial, because
one can find a metric $d$ such that  $X$
has infinite Hausdorff measure $H_1(X)$
 and lower limit in (\ref{i:main_rec}) is equal to $+\infty$.
In addition, it is easy to see that  the identity $H_1 (X) = 0$
can be replaced by stronger equality $H_{tg(t)} = 0$, where
$g(t)$ is some monotonically non--increasing function, $g(t) \to +\infty$ as $t\to 0+$.
\label{n:H_1}
\\

In Theorem \ref{t:general_case}
we use dense sets without arithmetic progression
to construct a dynamical system with slow time of multiple recurrence.
But the main idea of  Theorem \ref{t:F_corr} is the same.
Indeed,  let $A \subseteq {\bf N}$ be a set without arithmetic progressions and $D^* (A) > 0$.
%По теореме \ref{t:F_corr}
By Furstenberg's Correspondence Principle
there exists a dynamical system and a measurable set $E$, $\mu E >0$
such that for all natural numbers $n$ we have
$\mu (E \cap T^{-n} E \cap \dots \cap T^{-(k-1)n} E) = 0$.
In other words, we obtain the dynamical system without multiple recurrence.
(Certainly, this contradicts  Theorem \ref{t:F_deg}
and we can derive Szemer\'{e}di's Theorem from
the Theorem on multiple recurrence, see Statement  \ref{st:F_Sz_derive}).
Thus Theorem \ref{t:general_case} is a
quantitative analog of Theorem \ref{t:F_corr}.

In section $3$
%Besides that
we shall consider  a question concerning possible values of
one--dimensional recurrence constant $C(x)$.

The constructions which we use develop the approach of paper \cite{Bo}
and book \cite{Fu}.

\newpage

\refstepcounter{section}

{\bf \arabic{section}. Proof of Theorem \ref{t:general_case}.}

We need in the following simple Lemma.
In fact, this Lemma was proved in \cite{Tij}.
Another proof can be found in \cite{Tao_lect}.

\Lemma
{ \it
  Suppose $N$ is a natural number,
  $A$ is a nonempty subset of  $\Z_N$
  and $\_phi \ge 1$ is a real number.
  Then there exist residues  $a_1,\dots,a_l \in \Z_N$ and
  a
  %disjoint
  partition of ${\bf Z}_N$
  into sets
  $A_1, A_2, \dots, A_l$ and $B$ such  that \\
  $1)~$ $A_i \subseteq A+a_i$ for all  $i=1,\dots,l$. \\
  $2)~$ $|A_i| \ge |A| / \_phi$ for all  $i=1,\dots,l$. \\
  $3)~$ $|B| \le N / \_phi$.
%  \\
%  $4)~$ $ ( \bigsqcup_{i=1}^l A_i ) \bigsqcup B = {\bf Z}_N$.
}
\label{l:Tijdeman}
\\
\Proof
The proof of  the Lemma is a sort of an inductive process.
At the $n$th step of our process sets   $A_1,\dots, A_n$, residues   $a_1, \dots, a_n$
and auxiliary sets $B_1,\dots,B_n$ will be constructed.
Besides that we will have $B_1 \supseteq B_2 \supseteq \dots \supseteq B_n$.

Let $n=1$.
We set
$a_1 = 0$, $A_1 = A$ and $B_1 = \Z_N \setminus A_1$.
If \\
$|B_1| \le N / \_phi$ then Lemma \ref{l:Tijdeman} is proved.
Indeed, let $B=B_1$.
Clearly, the sets $A_1$, $B$ and the residues  $a_1$ satisfy  $1)$---$3)$.

Suppose at the $n$th  step of our procedure
the sets $A_1, \dots, A_n$ and residues  $a_1,\dots, a_n$ are constructed.
Let  $B_n = \Z_N \setminus (\bigsqcup_{i=1}^n A_i)$.
If $|B_n| \le N / \_phi$ then Lemma \ref{l:Tijdeman} is proved.
Indeed, put $B=B_n$.
Clearly, the sets $A_1,\dots, A_n$,  $B$ and residues $a_1,\dots,a_n$ satisfy $1)$---$3)$.

Let $|B_n| > N / \_phi$.
We have
\begin{equation}\label{f:tmp_5_13:13}
  \sum_{t\in \Z_N} |B_n \cap (A+t)| = |A| |B_n| \ge \frac{N |A|}{\_phi} \,.
\end{equation}
Hence there exists  $t \in \Z_N$ such that $|B_n \cap (A+t)| \ge |A| / \_phi$.
Put $a_{n+1} = t$ and $A_{n+1} = B_n \cap (A+a_{n+1})$.
Then for all  $i=1,\dots,n$ we get $A_{n+1} \cap A_i = \emptyset$.
Besides that  $A_{n+1} \subseteq A + a_{n+1}$.

For any  $i$ we have  $|A_i| \ge |A| / \_phi > 0$.
This implies that our process stops  at step $K$,  $K\le [N \_phi / |A|]+1$.
This completes the proof.

{\bf Proof of Theorem \ref{t:general_case}.}
Let $\a_m$ be an arbitrary non--increasing sequence of real numbers $\a_m \in (0,1)$,
$\a_m$ tends to zero as $m$ tends to infinity.
The function $\psi(n)$ is defined on positive integers.
Denote by the same letter $\psi(t)$
the result of linear extension of $\psi(n)$ to the entire ${\bf R}$.
We obtain the continuous monotonically increasing function.
Let $\_phi(t) = \sqrt{\psi(t)}$ and $\_phi^*(t) = \max \{ 1, \_phi(t) \}$.
Let also
$N_0 \le N_1 \le \dots \le N_m \le \dots~$
be a non--increasing sequence of integer numbers,
where $N_0 = 1$ and for all $m\ge 1$ we have
$
 N_m = \lceil \_phi^{-1} (2 \a_m^{-1} \_phi^*(2) N_0 N_1 \dots N_{m-1} ) \rceil .
$
Here $\_phi^{-1}$ is the inverse function.
We have $N_m \ge 2$, $m\ge 1$.

Let $X$ be a space of sequences $(x_1,x_2, \dots)$, $0\le x_i < N_i$, $i \ge 1$.
$C(a_1,\dots, a_l) = \{ x = (x_1, x_2, \dots) \in X ~:~ x_1 = a_1, \dots, x_l = a_l \}$
is called {\it elementary cylinder } of rank $l$.
We can associate with the sequence $x = (x_1,x_2, \dots) \in X$ the real number
$x \to \sum_{i=1}^{\infty} \frac{x_i}{N_0 N_1 \dots N_i} \in [0,1]$.
Thus, $X$ can be considered as the segment $[0,1]$.

Let $a$ be a positive integer, $N \in {\bf N}$.
By $a^{+} (N)$ define the number $a+1 \pmod{N}$.
Let $T$ be a transformation of the space $X$ into itself such that
$Tx = y$, $x = (x_1,x_2, \dots)$, $y = (y_1,y_2, \dots)$, where
$$
  y_1 = x_1^{+} (N_1)\,,
$$
\begin{displaymath}
  y_2 =
  \left\{ \begin{array}{ll}
    x_2^{+} (N_1),            &       \mbox{ if } x_1 + 1 = N_1
                              \\
    y_2,                      &       \mbox{ otherwise. }
  \end{array} \right.
\end{displaymath}
$$
\dots
$$
\begin{displaymath}
  y_m =
  \left\{ \begin{array}{ll}
    x_m^{+} (N_1),            &       \mbox{ if } x_1 + 1 = N_1, x_2 + 1 = N_2, \dots, x_{m-1} + 1 = N_{m-1}
                              \\
    y_m,                      &       \mbox{ otherwise. }
  \end{array} \right.
\end{displaymath}
$$
\dots
$$
The space $X$ has a natural group operation $+$.
We have $Tx = x+1$, where $1=(1,0,0,\dots)$.
Clearly,  $T$ preserves Haar measure $\mu$
and this Haar measure coincides with Lebesgue measure.
Elementary cylinder of rank $l$ has measure $1/ (N_0 N_1 \dots N_l)$.

Consider an arbitrary $N_s$.
By assumption there exists a non--empty set $A^{(N_s)} = A^{(s)} \subseteq {\bf Z}_{N_s}$
without arithmetic progressions of the length   $k$.
Using Lemma \ref{l:Tijdeman} for $A^{(s)}$ and $\_phi = \_phi(N_s)$,
we obtain the sets
%построим для любого $N_s$
$A_1^{(s)}, \dots, A_l^{(s)}$, $l=l(s)$
and $B^{(s)}$ satisfy  $1)$---$3)$.

Let $x,y \in X$, $x=(x_1, x_2,\dots )$, $y=(y_1,y_2,\dots)$.
Consider the function
$$
  d(x,y) = \{
                \psi^{-1} (N_0 \dots N_m) \rho(N_0 \dots N_m) \, \mbox{, where } m \mbox{ is the maximal integer s.t. }
$$
$$
                x_1 = y_1, \dots, x_{m-1} = y_{m-1}
                    \mbox{ and either there exists } i\in 1,2, \dots, l(m)
%                        \mbox{ such that }
                        \mbox{ s.t. }
$$
$$
                        x_m, y_m \in A_i^{(m)}
                    \mbox{ or }
                        x_m, y_m \in B^{(m)}
           \} \,,
$$
where $\psi^{-1} = 1/\psi$.
It is easy to see that $d(x,y)$ is a non--archimedean metric on $X$.
Let us consider Hausdorff measure $H_1$ on the space $X$.
Any elementary cylinder is a closed set therefore it is a Borel set
in the metric space $(X,d)$.
It follows that the measure $\mu$ is congruent to $H_1$.

Let us prove that $H_1(X) = 0$.
Let $\delta$ be an arbitrary positive number.
Since $\rho (N) \to 0$ as $N \to \infty$
it follows that
there exists a natural  $m$ such that
$
 \rho (N_0 \dots N_m) / \psi (N_0 \dots N_m) < \delta .
$
Let us consider the following partition of $X$ into
$$
  U_i (\v{a}) = \{ x = (x_1,x_2,\dots ) \in X ~:~ x_1 = a_1, \dots, x_{m-1} = a_{m-1}, x_m \in A_i^{(m)} \} \,,
      i=1,\dots, l(m)
$$
 and
$$
  B (\v{a}) = \{ x = (x_1,x_2,\dots ) \in X ~:~ x_1 = a_1, \dots, x_{m-1} = a_{m-1}, x_m \in B^{(m)} \} \,,
$$
where
$
 \v{a} \in [N_1] \m \dots \m [N_{m-1}] := F_{m-1}.
$
We obtain
$$
  X = \bigsqcup_{\v{a} \in F_{m-1}} \left( B (\v{a}) \bigsqcup \left( \bigsqcup_{i=1}^{l(m)} U_i (\v{a}) \right) \right)\,.
$$
For any $\v{a} \in F_{m-1}$ and any  $i\in 1,2,\dots, l(m)$ we have                     \\
$\diam U_i (\v{a}) \le \rho (N_0 \dots N_m) / \psi (N_0 \dots N_m) < \delta$.
Similarly, for any $\v{a} \in F_{m-1}$ we get
$
 \diam B (\v{a}) \le \rho (N_0 \dots N_m) / \psi (N_0 \dots N_m) < \delta
$.
Using $2)$ of Lemma \ref{l:Tijdeman}, we obtain
$l(m) \le N_m \_phi (N_m) / |A^{(m)}| = \_phi (N_m) / \rho (N_m)$.
Hence
$$
  H_1^{\delta} (X) \le
                        |F_{m-1}| \left( \frac{\_phi (N_m)}{\rho (N_m)} + 1 \right) \frac{\rho (N_0 \dots N_m)}{\psi (N_0 \dots N_m)}
                   \le
$$
$$
                   \le
                        2 N_0 \dots N_{m-1} \frac{\_phi (N_m) \rho (N_0 \dots N_m)}{\rho (N_m) \psi (N_0 \dots N_m)}
                   \le
                        2 N_0 \dots N_{m-1} \frac{\_phi (N_m)}{\psi (N_m)} \,.
$$
We have $N_m \ge \_phi^{-1} (2 \a_m^{-1} N_0 \dots N_{m-1})$.
This implies that
$
 2 N_0 \dots N_{m-1} \le \a_m \_phi (N_m).
$
Using this inequality, we get
$$
  H_1^{\delta} (X) \le \a_m \frac{\_phi^2 (N_m)}{\psi (N_m)} \le \a_m \,.
$$
Since $\a_m \to 0$ as $m \to \infty$ it follows that $H_1 (X) = 0$.

Let us prove (\ref{i:main_rec}).

Let
$$
 {\t{B}}^{(s)} = \{ x \in X ~:~ x_s \in B^{(s)} \} \subseteq X
       \quad     \mbox{ and }     \quad
 {\bf B} = \bigcap_{n=1}^{+\infty} \bigcup_{s\ge n}^{+\infty} {\t{B}}^{(s)} \,.
$$
We have
$ \mu (\t{B}^{(s)}) =  |B^{(s)}| / N_s  \le 1/ \_phi(N_s)$.
%$|B^{(s)}|  \le 1/ \_phi(N_s)$.
Since
$\_phi(N_s) \ge N_0 \dots N_{s-1}$, $s \ge 1$ and $N_s \ge 2$, $s\ge 1$
it follows that
\begin{equation}\label{}
  \sum_{s=1}^{\infty} \mu (\t{B}^{(s)})
    \le
  \sum_{s=1}^{\infty} \frac{1}{\_phi(N_s)} \le \sum_{s=1}^{\infty} \frac{1}{N_0 \dots N_{s-1}} < \infty \,.
\end{equation}
Using Borel--Cantelli Lemma, we get $\mu {\bf B} = 0$.

Let us prove that (\ref{i:main_rec}) holds for any $x \notin {\bf B}$.

Let $x = (x_1, x_2, \dots ) \in X \setminus {\bf B}$.
Since $x\notin {\bf B}$
it follows that
there exists a number $M = M(x) \in {\bf N}$ such that
for all $n\ge M$ we have $x_n \notin B^{(n)}$.
Let $m_0$ be a natural such that
$
 N_0 \dots N_{m_0-1} < M \le N_0 \dots N_{m_0}.
$

Prove that for any $n\ge N_0 \dots N_{m_0}$ the following inequality holds
\begin{equation}\label{i:tmp_ge}
  \rho^{-1} (n) \psi (n) \cdot \max  \{ d(T^n x,x), d(T^{2n} x,x), \dots, d(T^{(k-1)n} x,x) \} \ge 1 \,.
\end{equation}
Let $m_1 \ge m_0$ be a natural number.
Suppose that for some $n>0$ such that
\begin{equation}\label{c:ineq_n}
 N_0 \dots N_{m_1} \le n < N_0 \dots N_{m_1 + 1}
\end{equation}
(\ref{i:tmp_ge}) does not hold.
Then
$$
  d(T^n x, x),\, d(T^{2n} x, x),\, \dots,\, d(T^{(k-1)n} x, x) < \rho (N_0 \dots N_{m_1}) \psi^{-1} (N_0 \dots N_{m_1}) \,.
$$
Let $y^{(1)} = T^n x, y^{(2)} = T^{2n} x, \dots , y^{(k-1)} = T^{(k-1)n} x$.
Using properties of metric $d(x,y)$, we obtain
$$
  d(y^{(1)}, x), \dots, d(y^{(k-1)}, x) \le \rho (N_0 \dots N_{m_1+1}) \psi^{-1} (N_0 \dots N_{m_1+1}) \,.
$$
It follows that
$$
  x_1 = y^{(1)}_1 = \dots = y^{(k-1)}_1, \, \dots\,, x_{m_1} = y^{(1)}_{m_1} = \dots = y^{(k-1)}_{m_1}
$$
\begin{equation}\label{e:m_le}
            \mbox{ and there exists  }   i \mbox{ such that }
            x_{m_1+1}, y^{(1)}_{m_1+1}, \dots, y^{(k-1)}_{m_1+1} \in A^{(m_1+1)}_i \,.
\end{equation}
We have $n = y^{(1)}-x = y^{(2)}-y^{(1)} = \dots = y^{(k-1)} - y^{(k-2)}$.
Using (\ref{e:m_le}), we get
$$
  y^{(1)}-x = (\underbrace{0,\dots,0}_{m_1} \,, (y^{(1)}_{m_1+1} - x_{m_1+1}) \pmod{N_{m_1+1}}, w_1, w_2, \dots ) \,,
$$
$$
  y^{(2)}-y^{(1)} = (\underbrace{0,\dots,0}_{m_1} \,, (y^{(2)}_{m_1+1} - y^{(1)}_{m_1+1}) \pmod{N_{m_1+1}}, w'_1, w'_2, \dots ) \,,
  \dots \,,
$$
$$
  y^{(k-1)}-y^{(k-2)} = (\underbrace{0,\dots,0}_{m_1} \,, (y^{(k-1)}_{m_1+1} - y^{(k-2)}_{m_1+1}) \pmod{N_{m_1+1}}, w''_1, w''_2, \dots ) \,,
$$
where $w_1,w_2,\dots$, $w'_1, w'_2, \dots$ and $w''_1, w''_2, \dots$ are some numbers.
It follows that                                                                                     \\
$x_{m_1+1}, y^{(1)}_{m_1+1}, \dots, y^{(k-1)}_{m_1+1}$
is an arithmetic progression of length $k$
in ${\bf Z}_{N_{m_1+1}}$.
We have
$
 x_{m_1+1}, y^{(1)}_{m_1+1}, \dots, y^{(k-1)}_{m_1+1} \in A^{(m_1+1)}_i
$.
Since $A^{(m_1+1)}_i \subseteq A^{(m_1+1)}+p$ for some $p \in {\bf Z}_{N_{m_1+1}}$
it follows that
$A^{(m_1+1)}_i$ contains no arithmetic progressions of the length $k$
in ${\bf Z}_{N_{m_1+1}}$.
Hence for any $l=1,2, \dots, k-1$ we have
$x_{m_1+1} \equiv y^{(l)}_{m_1+1} \pmod{N_{m_1+1}}$.
Since $0\le x_{m_1+1}, y^{(1)}_{m_1+1}, \dots, y^{(k-1)}_{m_1+1} < N_{m_1+1}$
it yields that
$x_{m_1+1} =  y^{(1)}_{m_1+1} =  \dots = y^{(k-1)}_{m_1+1}$.
Hence
\begin{equation}\label{tmp_1:19:37}
  n = y^{(1)} - x = (\underbrace{0,\dots, 0}_{m_1+1} \,, n_{m_1+2}, n_{m_1+3}, \dots) \,.
\end{equation}
Using (\ref{tmp_1:19:37}), we get
$n \ge N_0 \dots N_{m_1+1}$.
This contradicts (\ref{c:ineq_n}).
Theorem \ref{t:general_case} is proved.

\Cor
{\it
 Let $k$ be an integer, $k\ge 3$.
 For any $\eps>0$ there exists an absolute positive constant $C_k$
 depends on  $k$ only
 and
 a dynamical system
 $(X,\mathcal{B},\mu,T,d)$,
 $X=[0,1]$, $\mu$ is Lebesgue measure
 such that
 $\mu$ is congruent to Hausdorff measure $H_1$,
 $H_1 (X) = 0$
 and
 for almost every point $x\in X$
 \begin{equation}\label{i:main_rec1}
    \liminf_{n\to \infty}~ \{ \rho^{-1} (n) \max  \{ d(T^n x,x), d(T^{2n} x,x), \dots, d(T^{(k-1)n} x,x) \} \} \ge 1 \,,
 \end{equation}
 where
 $
   \rho^{-1} (n) = 1 / \rho (n)
 $
 and $\rho(n) = \exp( - (1+\eps) C_k (\log n)^{1/(k+1)} )$.
}
\\
\Proof
By Theorem \ref{t:Be_Ra}
for any integer $k\ge 3$ and
for all sufficiently large $N$
there exists a set $A^{(N)}_0 \subseteq [1,2,\dots,N)$
contains no arithmetical progressions of the length  $k$ such that
$|A^{(N)}_0| \ge N \exp( -(1+\eps) C_k (\log N)^{1/(k+1)} )$, where
$C_k$ is an absolute positive constant depends on $k$ only.
Let
$A^{(N)}_1 = A^{(N)}_0 \cap [1,N/k)$,
$A^{(N)}_2 = A^{(N)}_0 \cap [N/k,2N/k), \dots,$
$A^{(N)}_{k} = A^{(N)}_0 \cap [N(k-1)/k ,N)$.
Any set $A^{(N)}_j$, $j\in [1,2,\dots,k]$
has no arithmetic progressions of the length  $k$ in $\Z_N$.
Clearly, there exists $j\in [1,2,\dots,k]$ such that
$|A^{(N)}_j| \ge \frac{N}{2k} \exp( -(1+\eps) C_k (\log N)^{1/(k+1)} )$.
Put $A^{(N)} = A^{(N)}_j$.
Using Theorem \ref{t:general_case}, we obtain the dynamical system
such that (\ref{i:main_rec}) holds for
$\rho(n) = \exp( -C_k (1+\eps') (\log n)^{1/(k+1)} )$, where $\eps'$
can be taken, for example, as $2\eps$.
This completes the proof.

\newpage

\refstepcounter{section}

{\bf \arabic{section}. On one--dimensional recurrence.}

\Th
{\it
 Let $f$ be a real number, $f\ge 1$,
 $X = [0,1]$ and $\mu$ be Lebesgue  measure on $X$.
 Then there exists a dynamical system
 $(X,\mathcal{B},\mu,T,d)$ such that
 $\mu$ is congruent to Hausdorff measure $H_1$,
 $H_1 (X) = 1$
 and
 for almost any point $x\in X$
 \begin{equation}\label{}
   C_f (x) := \liminf_{n\to \infty}~ \{ n \cdot f \cdot d(T^n x,x) \} = 1 \,.
 \end{equation}
}
\label{t:one_d_case}

\Note
Theorem \ref{t:one_d_case} was proved in \cite{Bo} in the case of
$f=1$.
\label{n:Bosh_f}

\Note
Let $X=[0,1]$, $\mathcal{B}$ be Borel $\sigma$--algebra, $\mu$ be Lebesgue measure,
$T_\a$ be a transformation of $X$ into itself, $T_\a x = (x+\a) \pmod{1}$.
Let also $d(x,y) = \| x-y \|$, where $\| \cdot \|$ is the integral distance.
There exists a number
$\lambda^*$,
$( \lambda^* = 5,68195..)$ such that for all
$f\ge \lambda^*$
there is  $T_\a$,
$\a = \a(f)$
such that $C_f (x) = 1$ for {\it all} $x \in [0,1]$.
The ray $[\lambda^*, +\infty)$ is called {\it Hall ray} (see \cite{Hall}).
The explicit value of  $\lambda^*$
was found by
G.R. Freiman
in
\cite{Freiman}.
By theorem \ref{t:one_d_case}
for all $f\ge 1$  (not only for $f\ge \lambda^*$)
there exists
a dynamical system with
$C_f (x) = 1$ for almost all $x \in [0,1]$.

\Note
The inequality  $H_1 (X) \le  1$ in Theorem \ref{t:one_d_case}
is important (see Note \ref{n:H_1}).
Besides that the inequality $H_1(X)\ge 1$ is very important too.
If this inequality does not hold then Theorem  \ref{t:one_d_case}
is trivial.
Indeed, let $f > 1$ and $(X,\mathcal{B},\mu,T,d)$
be the dynamical system such that $C_1 (x) = 1$ for all $x$
(see Note \ref{n:Bosh_f}).
Put $\t{d}(x,y) = d(x,y)/f$ and consider
the new dynamical system
$(X,\mathcal{B},\mu,T,\t{d})$.
Then for any $x\in X$ we have $C_f (x) = 1$.
Note that $H_1 (X) = 1/f < 1$ in this dynamical system.

\Proof
Let $N_0 = 1$,
$N_m = \lceil f 2^m \rceil^2$, $m = 1,2,\dots$
 and
let $X$ be a space of sequences
$(x_1,x_2, \dots)$, $0\le x_i < N_i$, $i \ge 1$.
There is a correspondence between $X$
and
$[0,1]$ is given by
$x \to \sum_{i=1}^{\infty} \frac{x_i}{N_0 N_1 \dots N_i}$.
Therefore we can consider the space $X$
to be the segment $[0,1]$.
Let also
the transformation
$T: X \to X$ is given by $Tx = x+1$,
where the addition was defined in the proof of Theorem \ref{t:general_case}
and $1=(1,0,0,\dots)$.
It was noted that $T$ preserves Lebesgue measure $\mu$.

Let $p_m = \sqrt{N_m} = \lceil f 2^m \rceil$, $m\ge 1$ and let                                              % Дописать
$$
  A_m^{(j)} = \{ x \in [0,1,\dots, N_m - 1] ~:~ x \equiv j \pmod{p_m} \} \,.
$$
Clearly, $[0,1,\dots, N_m - 1]$ is partitioned into the sets $A_m^{(j)}$, $j=0,1,\dots, p_m-1$.
Let the mapping
$
 \_phi_j : A_m^{(j)} \to {\bf N}_0 = {\bf N} \cup \{ 0 \}
$
are given by
$
 \_phi_j (x) = (x-j) / p_m.
$
If $x \in A_m^{(j)}$ then
we put
$
 \_phi (x) := \_phi_j (x).
$
It follows that the function $\_phi (x)$ is well--defined on $[0,1,\dots, N_m - 1]$.

Consider the function
$$
  d(x,y) = \Big{\{}
                \frac{ r_m (x_m,y_m) } {N_0 \dots N_{m-1}} \, \mbox{, where } m \mbox{ is the maximal integer such that }
$$
$$
                x_1 = y_1, \dots, x_{m-1} = y_{m-1}
                    \mbox{ and there exists } i\in 1,2, \dots, l(m)
                        \mbox{ such that }
$$
$$
                        x_m, y_m \in A_m^{(i)}
           \Big{\}} \,,
$$
where
\begin{displaymath}
  r_m (x_m, y_m)
                 =
  \left\{ \begin{array}{ll}
    \frac{1}{N_m},            &       \mbox{ if } x_m = y_m
                              \\
    \frac{ | \_phi(x_m) - \_phi(y_m) | }{f p_m},                      &       \mbox{ otherwise. }
  \end{array} \right.
\end{displaymath}
Note that
$
 1/N_m \le r_m (x_m,y_m) \le 1 .
$
\\
{\bf Statement.}
{\it $d(x,y)$ is a metric on $X$}.
\\
{\it Proof of the Statement.}
Clearly, $d(x,y)$ is a symmetric function
and  $d(x,y) = 0$ if and only if  $x=y$.
Let us prove that for any  $x,y,z \in X$,
$x = (x_1,x_2, \dots )$, $y = (y_1,y_2, \dots )$, $z = (z_1,z_2, \dots )$
we have
\begin{equation}\label{i:tr}
  d(x,y) \le d(x,z) + d(z,y) \,.
\end{equation}
If $d(x,y) = 0$ then
(\ref{i:tr}) is trivial.
Suppose $d(x,y) > 0$.
Then there exists $m\in {\bf N}$ such that
$d(x,y) = r_m (x_m, y_m) / (N_0 \dots N_{m-1})$.
If there exists $i \in 1,2, \dots, m-1$ such that either
$z_i \neq x_i$ or $z_i \neq y_i$ then
(\ref{i:tr}) holds.
Therefore we can assume that for all   $i=1,2, \dots, m-1$ we have $z_i = x_i = y_i$.

Suppose that for any  $j$ the elements
$x_m, y_m$ do not belong to
$A_m^{(j)}$.
It follows that $d(x,y) = 1/ (N_0 \dots N_{m-1})$.
On the other hand, either $x_m, z_m$ or $y_m, z_m$ do not belong to the same set $A_m^{(j)}$.
Hence (\ref{i:tr}) holds again.

Suppose  $x_m$ and $y_m$ belong to the same set $A_m^{(j)}$.
If $z_m \notin A_m^{(j)}$ then (\ref{i:tr}) holds.
If $z_m \in A_m^{(j)}$ then
$$
  d(x,y) = \frac{| \_phi(x_m) - \_phi(y_m) |}{N_0 \dots N_{m-1} f p_m}
           \le
           \frac{| \_phi(x_m) - \_phi(z_m) |}{N_0 \dots N_{m-1} f p_m}
           +
           \frac{| \_phi(z_m) - \_phi(y_m) |}{N_0 \dots N_{m-1} f p_m}
           =
$$
$$
           =
           d(x,z) + d(z,y) \,.
$$
This completes the proof of the Statement.

Let us return to the proof of Theorem  \ref{t:one_d_case}.
Consider Hausdorff measure $H_1$ on $X$.
Any elementary cylinder is a closed set therefore it is a Borel set
in the metric space $(X,d)$.
Hence the measure $\mu$ is congruent to $H_1$.
We claim that $H_1(X) = 1$.

Consider the sets
$$
  C(a_1,\dots, a_m) = \{ x = (x_1,x_2, \dots ) \in X ~:~ x_1 = a_1, \dots, x_m = a_m \} \,.
$$
Obviously, the space $X$ is partitioned into these sets.
Hence $H_1(X) \le 1$.

Let us prove that  $H_1(X) \ge 1$.
Suppose $H_1 (X) = a < 1$.
Since
$H_1 (X) = \lim_{\d \to 0} H_1^{\d} (X) = \sup_{\d>0} H_1^{\d} (X)$
(see \cite{Falconer} for example) it follows that
for any $\eps>0$
there exists $\d>0$ such that
\begin{equation}\label{tmp:15:23:26}
  a-\eps < H_1^{\d} (X) \le a = H_1 (X) \,.
\end{equation}
Let $\eps_0 = (1-a)/2 > 0$.
Using (\ref{tmp:15:23:26}) and the definition of Hausdorff measure,
we obtain a covering of $X$ by sets  $\{ U_i \}$, $r_i = \diam U_i$,
$r_i < \d = \d(\eps_0)$ such that
\begin{equation}\label{tmp:15:23:29}
  a-\eps < \sum_i \diam U_i = \sum_i r_i < a+ \eps \,.
\end{equation}
If $a=0$ then the left--hand side  of (\ref{tmp:15:23:29})
is not really need.

If $r_i = 0$ then the set $U_i$ is a one--point set,
$U_i = \{ p_i \}$.
Denote by $P$ the union of all one--point sets $U_i$.
In other words, $P = \cup_{ \{i: r_i = 0\} } U_i = \cup_i \{ p_i \}$.
Clearly, there exists  $U_i$ does not belong to $P$.
We shall consider only these sets $U_i$.
Since zero is a unique limit point of
the set of distances of $X$
it follows that for any $U_i$
there exist two points
$x,y \in U_i$ such that
$r_i = \diam U_i = d(x,y)$.
Let $d(x,y) = r_m (x_m, y_m) / (N_0 \dots N_{m-1})$.
If there exists a number  $j$, $j=j(i)$ such that $x_m,y_m \in A_m^{(j)}$
then put
\begin{equation}\label{e:tmp:1}
  C_i = \{ z = (z_1,z_2,\dots) \in X ~:~ z_1 = x_1, \dots, z_{m-1} = x_{m-1}, z_m \in  A_m^{(j)} \cap [x_m, y_m] \} \,.
\end{equation}
If there is not such $A_m^{(j)}$ then we set
\begin{equation}\label{e:tmp:2}
  C_i = \{ z = (z_1,z_2,\dots) \in X ~:~ z_1 = x_1, \dots, z_{m-1} = x_{m-1} \} \,,
\end{equation}
Clearly, in
%the
both cases we have
$U_i \subseteq C_i$ and $\diam C_i = \diam U_i$.
It follows that $C_i$ satisfy (\ref{tmp:15:23:29})
and  $\{ \{ C_i \}, P \}$ is a partition of  $X$.

Note that if $C_i$ is given by (\ref{e:tmp:2}) then
$C_i$ is an elementary cylinder.
Let
$$
  C_i (a) = \{ z = (z_1,z_2,\dots) \in X ~:~ z_1 = x_1, \dots, z_{m-1} = x_{m-1}, z_m = a \} \,.
$$
Then $C_i (a)$ is an elementary cylinder for all $a$, $0\le a < N_m$.
If $C_i$ is given by (\ref{e:tmp:2}) then
$C_i = \bigsqcup_{a\in A_m^{(j)} \cap [x_m, y_m]} C_i (a)$.
Clearly,
$\diam C_i \ge \sum_{a\in A_m^{(j)} \cap [x_m, y_m]} \diam C_i (a)$.
It follows that there exists a countable covering of $X$ by $P$
and elementary cylinders $C'_i$, $r'_i = \diam C'_i$
such that
\begin{equation}\label{}
  \sum_i r'_i \le \sum_i r_i < a + \eps < 1 \,.
\end{equation}
Suppose $C_1$ and $C_2$ are two elementary cylinders.
Then either $C_1$, $C_2$  are disjoint or one of them contains another.
Therefore there exists a sub--covering $C''_i$, $r''_i = \diam C''_i$ of the covering $C'_i$
such that $C''_i$ and $P$ is a {\it partition} of $X$ by elementary cylinders
and
\begin{equation}\label{tmp:13:46}
 \sum_i r''_i \le \sum_i r'_i < 1 \,.
\end{equation}
We have  $r''_i = \mu C''_i$ and
$\sum_i r''_i = \sum_i \mu C''_i = \mu(X) = 1$.
This contradicts  (\ref{tmp:13:46}).
Hence $H_1 (X) = 1$.

Finally, we need to prove for almost all $x\in X$ the following inequality
\begin{equation}\label{e:desire1}
   \liminf_{n\to \infty}~ \{ n \cdot f \cdot d(T^n x,x) \} = 1 \,.
\end{equation}
By $a_m (j)$ denote the maximal element of $A_m^{(j)}$.
Let $B_m = \bigsqcup_{j=1}^{p_m} a_m (j)$.
We have
$|B_m| = p_m = \sqrt{N_m}$.
Let
$$
 \t{B}_m = \{ x \in X ~:~ x_m \in B_m \} \subseteq X
       \quad     \mbox{ and }     \quad
 {\bf B} = \bigcap_{n=1}^{+\infty} \bigcup_{m\ge n}^{+\infty} \t{B}_m \,.
$$
Then
$ \mu (\t{B}_m) =  |B_m| / N_m =  1/ p_m$.
We have
\begin{equation}\label{}
  \sum_{m=1}^{\infty} \mu (\t{B}_m) =
  \sum_{m=1}^{\infty} \frac{1}{\sqrt{N_m}} \le \sum_{m=1}^{\infty} \frac{1}{2^m} < \infty \,,
\end{equation}
Using Borel--Cantelli Lemma, we get $\mu {\bf B} = 0$.

Let us prove (\ref{e:desire1}) for all
$x=(x_1,x_2,\dots)$, $x\notin {\bf B}$.
If $x \notin {\bf B}$ then there exists $M = M(x)$
such that for all $n \ge M$ we have $x_n \notin B_n$.
We can assume that $M$ is a sufficiently large number.
There exists a natural $m_0$ such that $N_0 \dots N_{m_0} \ge M$.
Consider the increasing sequence of natural numbers
$$
  S = \{ p_{m+1} N_0 \dots N_m \}_{m=m_0}^{+\infty} = \{ n_m \}_{m=m_0}^{+\infty} \,.
$$
Let $x=(x_1,x_2, \dots, x_m, x_{m+1}, x_{m+2}, \dots)$, where $x_{m+1}$
belongs to some $A_{m+1}^{(j)}$.
Let also $n_m \in S$.
Then
$T^{n_m} x = (x_1,\dots, x_m, \t{x}_{m+1}, \t{x}_{m+2}, \dots)$,
where
$\t{x}_{m+1}, \t{x}_{m+2}, \dots$ are numbers such that
$x_{m+1}, \t{x}_{m+1}$ belong to $A_{m+1}^{(j)}$ and
$|\_phi(x_{m+1}) - \_phi(\t{x}_{m+1})| = 1$.
It follows that $d(T^{n_m} x,x) = 1/(N_0 \dots N_m f p_{m+1})$.
Further,
\begin{equation}\label{}
  n_m f \cdot d(T^{n_m} x,x) = p_{m+1} N_0 \dots N_m f \frac{1}{N_0 \dots N_m f p_{m+1}} = 1 \,.
\end{equation}
Hence for all $x \notin {\bf B}$ we have $C_f (x) \le 1$.

Prove that for all $x \notin {\bf B}$ the inverse inequality holds :  $C_f (x) \ge 1$.
Let $n$ be a natural number such that
$n \in [N_0 \dots N_m, N_0 \dots N_{m+1}) := J_m$ and $m \ge m_0$.
Note that $n_m$ belongs to  $J_m$.
If $n = t N_0 \dots N_m$, $1\le t < N_{m+1}$ then
$T^{n} x = (x_1,\dots, x_m, \t{x}_{m+1}, \t{x}_{m+2}, \dots)$,
where
$\t{x}_{m+1}, \t{x}_{m+2}, \dots$ are numbers.
Suppose  $\t{x}_{m+1} \notin A_{m+1}^{(j)}$
it follows that
$
 d(T^n x,x) = 1/(N_0 \dots N_m)
$
and consequently,
$
 n f d(T^n x,x) \ge 1.
$
If
$\t{x}_{m+1} \in A_{m+1}^{(j)}$ then
$1 = n_m f d(T^{n_m} x, x) \le n f d(T^n x, x)$.

Finally, suppose that
$n \neq t N_0 \dots N_m$, $1\le t < N_{m+1}$.
In the case
$T^{n} x = (x'_1,\dots, x'_m, x'_{m+1}, x'_{m+2}, \dots)$
and there exists $i\in 1,2, \dots, m$ such that $x_i \neq x'_i$.
It follows that
$
 d(T^n x,x) \ge 1/(N_0 \dots N_m)
$
 and again
$
 n f d(T^n x,x) \ge 1.
$

Thus for any $m\ge m_0$ and an arbitrary $n\in [N_0 \dots N_m, N_0 \dots N_{m+1})$
we have
$
 1 \le n f d(T^n x,x) .
$
Whence for any $x\notin {\bf B}$
%the following inequality holds
we obtain
$C_f (x) \ge 1$.
This completes the proof.

\end{document}